\documentclass{amsart}
\usepackage{amssymb}
\usepackage{amsmath}

\usepackage{diagrams}
\diagramstyle[centredisplay,dpi=600,UglyObsolete,heads=LaTeX]

\newcommand\Zset{\mathbb {Z}}
\newcommand\Qmax{Q^r_{\mathrm{max}}}
\newcommand\Qlmax{Q^l_{\mathrm{max}}}
\newcommand\Qsimmax{Q^{\sigma}_{\mathrm{max}}}
\newcommand\Qtot{Q^r_{\mathrm{tot}}}

\newcommand\Qsimtot{Q^{\sigma}_{\mathrm{tot}}}
\newcommand\Qgold{Q^r_{\mathrm{Gol}}}
\newcommand\Qsimgold{Q^{\sigma}_{\mathrm{Gol}}}
\newcommand\Qcl{Q^r_{\mathrm{cl}}}
\newcommand\Qlrcl{Q_{\mathrm{cl}}}
\newcommand\f{{\mathcal F}}
\newcommand\te{{\mathcal T}}
\newcommand\ef{{\mathfrak F}}
\newcommand\de{{\mathfrak D}}
\newcommand\coker{{\mathrm{coker}}}
\newcommand\dirlim{\mathop{\varinjlim}\limits}
\newcommand\homo{\mathrm{Hom}}
\newcommand\ann{{\mathrm{ann}}}

\newtheorem{theorem}{Theorem}
\newtheorem{lemma}{Lemma}
\newtheorem{corollary}{Corollary}
\newtheorem{prop}{Proposition}
\newtheorem{definition}{Definition}
\newtheorem{example}{Example}

\begin{document}

\title{EXTENDING RING DERIVATIONS TO RIGHT AND SYMMETRIC RINGS AND MODULES OF QUOTIENTS}

\author{\footnotesize LIA VA\v{S}}

\address{Department of Mathematics, Physics and Statistics\\ University of the Sciences
in Philadelphia\\ Philadelphia, PA 19104, USA}

\email{l.vas@usip.edu}



\maketitle


\begin{abstract}
We define and study the symmetric version of differential torsion theories. We prove that the symmetric versions of some of the existing results on derivations on right modules of quotients hold for derivations on symmetric modules of quotients.
In particular, we prove that the symmetric Lambek, Goldie and perfect torsion theories are differential.

We also study conditions under which a derivation on a right or symmetric module of quotients extends to a right or symmetric module of quotients with respect to a larger torsion theory. Using these results, we study extensions of ring derivations to maximal, total and perfect right and symmetric rings of quotients.
\end{abstract}

\section{Introduction}

Recall that a {\em derivation} on a ring $R$ is a mapping $\delta: R \rightarrow R$ such that $\delta(r+s)=\delta(r)+\delta(s)$ and $\delta(rs)=\delta(r)s+r\delta(s)$ for all $r,s\in R.$ A mapping $d: M\rightarrow M$ on a right $R$-module $M$ is a {\em $\delta$-derivation} if $d(x+y)=d(x)+d(y)$ and $d(xr)=d(x)r+x\delta(r)$ for all $x\in M$ and $r\in R.$  In \cite{Golan_paper}, \cite{Bland_paper}, and \cite{Lia_Diff}, the authors study how derivations agree with an arbitrary hereditary torsion theory for that ring. In this paper, we continue that study.

Through this paper, we shall use the usual definition of torsion theory and hereditary torsion theory (e.g. \cite{Stenstrom}, \cite{Bland_book}, \cite{Bland_paper}, \cite{Lia_Diff}). If $\tau=(\te, \f)$ is a torsion theory for $R$ with $\te$ the class of torsion modules and $\f$ the class of torsion-free modules, we denote the
torsion submodule of a right $R$-module $M$ by $\te(M)$ and the torsion-free quotient $M/\te(M)$ by $\f(M).$ If $\tau$ is hereditary, we denote its Gabriel filter by $\ef.$ If $\te(R)=0,$ $\tau$ is said to be faithful.

A Gabriel filter $\ef$ is a {\em differential filter} if for every $I\in \ef$ there is $J\in \ef$ such that $\delta(J)\subseteq I$ for all derivations $\delta.$ The hereditary torsion theory determined by $\ef$ is said to be {\em differential} in this case. By Lemma 1.5 from \cite{Bland_paper}, a torsion theory is differential if and only if
\[d(\te(M))\subseteq \te(M)\]
for every right $R$-module $M,$ every derivation $\delta$ and every $\delta$-derivation $d$ on $M.$

If $\tau$ is a hereditary torsion theory with Gabriel filter $\ef$ and $M$ is a right $R$-module, the module of quotients $M_{\ef}$ of $M$ is defined as the largest submodule $N$ of the injective envelope $E(M/\te(M))$ of $M/\te(M)$ such that $N/(M/\te(M))$ is torsion module (i.e. the closure
of $M/\te(M)$ in $E(M/\te(M))$). Equivalently, the module of quotients $M_{\ef}$ can be defined by
\[M_{\ef}=\dirlim_{I\in\ef}\homo(I, \frac{M}{\te(M)})\]
(see chapter IX in \cite{Stenstrom}). Note that from this description it directly follows that $M_{\ef}=(M/\te(M))_{\ef}.$

The $R$-module $R_{\ef}$ has a ring structure and $M_{\ef}$ has a structure of a right
$R_{\ef}$-module (see exposition on pages 195--197 in \cite{Stenstrom}). The ring $R_{\ef}$ is called the
right ring of quotients with respect to the torsion theory $\tau.$

Consider the map $q_M:M\rightarrow M_{\ef}$ obtained by
composing the projection $M\rightarrow M/\te(M)$ with the injection $M/\te(M)\rightarrow M_{\ef}.$ The kernel and cokernel of $q_M$ are torsion modules and $M_{\ef}$ is torsion-free (Lemmas 1.2 and 1.5, page 196, in \cite{Stenstrom}). The maps $q_M$ define a left exact functor $q$ from the category of right $R$-modules to the category of right $R_{\ef}$-modules (see \cite{Stenstrom} pages 197--199). The functor $q$ maps torsion modules to 0 (see Lemma 1.3, page 196 in  \cite{Stenstrom}).

In Theorem on page 277 and Corollary 1 on page 279 of \cite{Golan_paper}, Golan has shown the following
\begin{prop}[Golan] Let $\delta$ be a derivation on $R$, $M$ a right $R$-module, $d$ a $\delta$-derivation on $M$ and $\tau=(\te, \f)$ a hereditary torsion theory.
\begin{enumerate}
\item If $M$ is torsion-free, then $d$ extends to a derivation on the module of quotients $M_{\ef}$ such that  $d q_M=q_M d.$

\item If $d(\te(M))\subseteq \te(M),$ then $d$ extends to a derivation on the module of quotients $M_{\ef}$ such that $d q_M=q_M d.$
\end{enumerate}
\label{Golan_Proposition}
\end{prop}

A direct corollary of the first part of Proposition \ref{Golan_Proposition} is that a ring derivation extends to a right ring of quotients for every hereditary and faithful torsion theory. By the second part of Proposition \ref{Golan_Proposition}, we can extend a derivation on a module to a derivation on its module of quotients for every differential torsion theory. Bland proved that such extension is unique and that the converse is also true. Namely, Propositions 2.1 and 2.3 of his paper \cite{Bland_paper} state the following.

\begin{theorem}[Bland] Let $\ef$ be a Gabriel filter.
\begin{enumerate}
\item If a derivation on a module $M$ extends to a derivation on the module of quotients $M_{\ef},$ then such extension is unique.

\item The filter $\ef$ is differential if and only if
every derivation on any module $M$ extends uniquely to a derivation on the module of quotients $M_{\ef}.$
\end{enumerate}
\label{Bland_Theorem}
\end{theorem}

In this paper, first we shall address the following question: assuming that a ring derivation can be extended to a ring of quotients $Q_1$ and that $Q_1$ is contained in another ring of quotients $Q_2,$ can a derivation on $Q_1$ be extended to $Q_2?$ More generally, when does a derivation on a module of quotients extends to a module of quotients with respect to a larger torsion theory? We address these questions in section \ref{section_extending_on_right}. We list some of the conditions under which the extensions described above are possible (Proposition \ref{Extending_to_right_modules_of_quotients} and Corollary \ref{Extending_to_right_rings_of_quotients}). Then we study the extensions of ring derivations to maximal, classical, total and perfect right rings of quotients (Corollary \ref{extending_on_Qtot_and_Qmax}).

The maximal symmetric rings of quotients emerged (first introduced by Utumi in \cite{Utumi}, studied in \cite{Lanning} and \cite{Ortega_paper_Qsimmax}) as an attempt to introduce the symmetric version of the maximal right (and left) ring of quotients. Schelter's work on symmetric rings of quotients in \cite{Schelter} parallels the work on a right ring of quotients with respect to a torsion theory -- it provides the basis for a uniform treatment of two-sided rings of quotients using torsion theories. Namely, if $\ef_l$ and $\ef_r$ are left and right Gabriel filters, the symmetric filter $_{\ef_l}\ef_{\ef_r}$ and the symmetric ring of quotients $_{\ef_l}R_{\ef_r}$ can be defined. In \cite{Ortega_thesis} and \cite{Ortega_paper}, Ortega defines the symmetric module of quotients $_{\ef_l}M_{\ef_r}$ of an $R$-bimodule $M.$

In section \ref{section_differentiability_of_symmetric}, we prove that the symmetric versions of Golan's and Bland's results on derivations on right modules of quotients hold for derivations on symmetric modules of quotients (Proposition \ref{symmetric_differential_TT} and Theorem \ref{extending_symmetric_diff})

In \cite{Lia_Sym}, the symmetric version of a right perfect ring of quotients and the symmetric version of total right ring of quotients are defined and studied. In this paper, we study their differentiability and prove the result analogous to Proposition 4 of \cite{Lia_Diff} -- we prove that a perfect symmetric filter is differential (Corollary \ref{symmetric_perfect}).

In \cite{Lia_Diff}, it has been proven that the Lambek and Goldie torsion theories for every ring are differential (Proposition 9 and 14 of \cite{Lia_Diff}). In this paper, we prove that the same is true for the symmetric version of Lambek and Goldie torsion theories (Corollary \ref{symmetric_Lambek_Goldie}).

In section \ref{section_extending_on_symmetric}, we prove the symmetric version of results on right torsion theories from section \ref{section_extending_on_right} (Proposition \ref{Extending_to_right_modules_of_quotients} and Corollaries \ref{Extending_to_right_rings_of_quotients} and \ref{extending_on_Qtot_and_Qmax}) -- we study the conditions under which a derivation on symmetric module of quotients extends to the symmetric module of quotients with respect to a larger symmetric torsion theory (Proposition \ref{extending_on_symmetric_modules}). Using this result, we study the extensions of ring derivations to maximal, total and perfect symmetric rings of quotients (Corollary \ref{extending_on_symmetric_Lambek_Goldie_perfect}).

\section{Extending derivations to right rings and modules of quotients}
\label{section_extending_on_right}

Let us suppose that $Q_1$ and $Q_2$ are two rings of quotients of a ring $R$ with derivation $\delta$ such that $Q_1\subseteq Q_2.$ The two questions are of interest.
\begin{itemize}
\item[Q1] If $\delta$ extends to $Q_1,$ can we extend it to $Q_2$ as well?

\item[Q2] If $\delta$ extends to both $Q_1$ and $Q_2,$ is the extension on $Q_2$ equal to the extension on $Q_1$ when restricted to $Q_1$? In other words, does the following diagram commute?
\end{itemize}

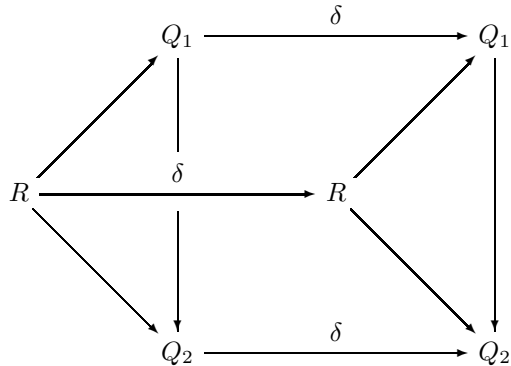
\begin{diagram}
  &        & Q_1           &     & \rTo^{\delta} &       & Q_1\\
  & \ruTo  & \vLine        &     &               & \ruTo &    \\
R &        & \rTo^{\delta} &     & R             &       & \dTo\\
  & \rdTo  & \dTo          &     &               & \rdTo &     \\
  &        & Q_2           &     & \rTo^{\delta} &       & Q_2\\
\end{diagram}

If the diagram above commutes, we shall say that the extensions on $Q_1$ and $Q_2$ {\em agree.}

More generally, let us consider modules of quotients. Suppose that a Gabriel filter $\ef_1$ is contained in a Gabriel filter $\ef_2$ and that $M$ is a right $R$-module. Let $\tau_i$ denote the corresponding torsion theories and $q_i$ denote the left exact functors mapping $M$ to the rings of quotients $M_{\ef_i}$ with respect to $\ef_i$ for $i=1$ and $2.$ Recall that the maps $q_i$ are defined by $m\mapsto L_{m+\te_i(M)},$ $i=1,2$ where $L_{m+\te_i(M)}$ denotes the left multiplication with $m+\te_i(M).$

Note that we have a mapping $q_{12}: M_{\ef_1}\rightarrow M_{\ef_2}$ given by the composition of the natural homomorphisms
\[
\dirlim_{I\in \ef_1} \homo(I, \frac{M}{\te_1(M)})\rightarrow\dirlim_{I\in \ef_2} \homo(I, \frac{M}{\te_1(M)})\rightarrow\dirlim_{I\in \ef_2} \homo(I, \frac{M}{\te_2(M)})
\] since $\ef_1\subseteq \ef_2.$
Note also that $q_{12}q_1=q_2$ since $q_{12}(q_1(m))=q_{12}(L_{m+\te_1(M)})=L_{m+\te_2(M)}=q_2(m).$

Applying $q_2$ to the diagram
\begin{diagram}
0 & \rTo & \te_1(M) & \rTo & M & \rTo^{q_1} & M_{\ef_1}& \rTo & \coker q_1\\
& & \dTo_{\subseteq} & & \dTo_{=} & & \dTo_{q_{12}} & &\dTo\\
0 & \rTo & \te_2(M) & \rTo & M & \rTo^{q_2} & M_{\ef_2}& \rTo & \coker q_2
\end{diagram}
we obtain the diagram
\begin{diagram}
0 & \rTo & 0 & \rTo & M_{\ef_2} & \rTo & (M_{\ef_1})_{\ef_2}& \rTo &0\\
& & \dTo & & \dTo_{=} & & \dTo & & \dTo\\
0 & \rTo & 0 & \rTo & M_{\ef_2} & \rTo & (M_{\ef_2})_{\ef_2} & \rTo&  0
\end{diagram}
Since the functor $q_2$ maps all $\tau_2$-torsion modules to 0 and $\te_1\subseteq\te_2,$ the modules $\te_i(M)$ and $\coker q_i$ for $i=1,2$ are mapped to 0 by $q_2.$ Thus, by the commutativity of the diagram above, $(M_{\ef_1})_{\ef_2}$ is isomorphic to $M_{\ef_2}.$
The composition of $q_2: M_{\ef_1} \rightarrow (M_{\ef_1})_{\ef_2}$ with isomorphism $i: (M_{\ef_1})_{\ef_2}\rightarrow M_{\ef_2}$ satisfies $iq_2=q_{12}.$

\begin{definition}
Suppose that there is a $\delta$-derivation $d$ defined on a right $R$-module $M.$ If $d$ extends to derivations $d_1$ on $M_{\ef_1}$ and $d_2$ on $M_{\ef_2}$ such that the diagram
\begin{diagram}
  &             & M_{\ef_1}    &      & \rTo^{d_1} &             & M_{\ef_1}\\
  & \ruTo^{q_1} & \vLine       &      &            & \ruTo^{q_1} &    \\
M &             & \rTo^{d}     &      & M          &             & \dTo_{q_{12}}\\
  & \rdTo^{q_2} & \dTo_{q_{12}}&      &            & \rdTo^{q_2} &     \\
  &             & M_{\ef_2}    &      & \rTo^{d_2} &             & M_{\ef_2}\\
\end{diagram}
commutes, we say that the extensions of $d$ on $M_{\ef_1}$ and $M_{\ef_2}$ {\em agree.}
\label{def_of_agreement}
\end{definition}

\begin{lemma} Suppose that a Gabriel filter $\ef_1$ is contained in a Gabriel filter $\ef_2$ and that $M$ is a right $R$-module with a $\delta$-derivation $d$. If $d$ can be extended to $M_{\ef_1}$ and either
\begin{itemize}
\item[i)] $d$ can be extended from $M_{\ef_1}$ to $M_{\ef_2},$  or

\item[ii)] $d$ can be extended from $M$ to $M_{\ef_2},$
\end{itemize}
then the extensions of $d$ to $M_{\ef_1}$ and $M_{\ef_2}$ agree.
\label{lemma_for_extending_right}
\end{lemma}
\begin{proof}
Since $\ef_1\subseteq\ef_2,$ $q_{12}q_1=q_2.$ By assumption that $d$ can be extended to a derivation $d_1$ on $M_{\ef_1},$ $d_1q_1=q_1 d.$

In case i), we have that $d_2q_{12}=q_{12} d_1$ and need to prove that $d_2q_2=q_2 d.$ This is the case because
\[d_2q_2=d_2 q_{12} q_1 = q_{12} d_1 q_1 = q_{12} q_1 d = q_2 d.\]

In case ii), we have that $d_2q_2=q_2 d$ and need to prove that $d_2q_{12}=q_{12} d_1.$ Let $q\in M_{\ef_1}.$ Then there is a right ideal $I$ in $\ef_1$ such that $qI\subseteq M/\te_1(M).$ In this case $q_{12}(q)I\subseteq M/\te_2(M)$ by definition of map $q_{12}.$ Note also that if $qi=m+\te_1(M)$ for some $i\in I$ and $m\in M,$ then $q_{12}(q)i=m+\te_2(M)$ since
\[L_{q_{12}(q)i}=q_{12}(L_{qi})=q_{12}(L_{m+\te_1(M)})=q_{12}(q_1(m))=q_2(m)=L_{m+\te_2(M)}.\]

Note that the fact that derivation on $M$ can be extended to $M_{\ef_1}$ gives us that $d(\te_1(M))\subseteq\te_1(M)$ by part 2 of Theorem \ref{Bland_Theorem}. Thus $d$ defines a derivation $\overline{d_1}$ on $M/\te_1(M)$ such that $\overline{d_1}(m+\te_1(M))=d(m)+\te_1(M).$ The extension $d_1$ coincides with the extension of the derivation $\overline{d_1}$ from $M/\te_1(M)$ to $(M/\te_1(M))_{\ef_1}=M_{\ef_1}$ (see Corollary 1 in \cite{Golan_paper}). Thus $q_1(d(m))=q_1(\overline{d_1}(m+\te_1(M)))=d_1(q_1(m)).$
Similarly, $d$ defines a derivation on $M/\te_2(M)$  we shall call $\overline{d_2}$ with $q_2(d(m))=q_2(\overline{d_2}(m+\te_2(M)))=d_2(q_2(m)).$

Thus, the extension $d_1$ on $M_{\ef_1}$ is defined such that $d_1(q)i=d_1(qi)- q\delta(i) =L_{\overline{d_1}(qi)}-q\delta(i)$ for all $i\in I.$
As $I$ is in $\ef_2$ also and $q_{12}(q)I\subseteq M/\te_2(M)$, we have that $d_2(q_{12}(q))i= d_2(q_{12}(q)i)-q_{12}(q)\delta(i)=L_{\overline{d_2}(q_{12}(q)i)}-q_{12}(q)\delta(i).$

We show first that $d_2 (q_{12}(q))i=q_{12}(d_1(q))i.$
\[
\begin{array}{rcll}
d_2 (q_{12}(q))i & = & L_{\overline{d_2}(q_{12}(q)i)}-q_{12}(q)\delta(i) & (\mbox{by remark above})\\
& =&  L_{\overline{d_2}(m+\te_2(M))}-q_{12}(q)\delta(i) & (q_{12}(q)i=m+\te_2(M)
)\\
& =&  d_2(L_{m+\te_2(M)})-q_{12}(q)\delta(i) & (\mbox{see above}
)\\
& =&  d_2(q_2(m))-q_{12}(q)\delta(i) & (\mbox{by definition of }q_2)\\
& =& q_2(d(m))-q_{12}(q)\delta(i) & (d_2q_2 =q_2 d)\\
& =& q_{12}(q_1(d(m)))-q_{12}(q)\delta(i) & (q_2=q_{12}q_1)\\
& =&  q_{12}(d_1(q_1(m)))-q_{12}(q)\delta(i) & (d_1q_1 =q_1 d)\\
& =&  q_{12}(d_1(L_{m+\te_1(m)}))-q_{12}(q)\delta(i) & (\mbox{by definition of }q_1)\\
& =&  q_{12}(L_{\overline{d_1}(m+\te_1(m))})-q_{12}(q)\delta(i) & (\mbox{see above}
)\\
& =&  q_{12}(L_{\overline{d_1}(qi)})-q_{12}(q)\delta(i) & (qi=m+\te_1(M))\\
& =&  q_{12}(L_{\overline{d_1}(qi)}-q\delta(i)) & (q_{12}\mbox{ is an }R\mbox{-map})\\
& =&  q_{12}(d_1(q)i) & (d_1(q)i=L_{\overline{d_1}(qi)}-q\delta(i))\\
& =&  q_{12}(d_1(q))i & (q_{12}\mbox{ is an }R\mbox{-map})\\
\end{array}
\]

Thus, the left multiplication with $d_2 (q_{12}(q))-q_{12}(d_1(q))$ defines the zero $R$-map $I\rightarrow M_{\ef_2}.$ As $I\in\ef_2,$ this map extends to a $R$ map $f: R\rightarrow M_{\ef_2}$ (see Proposition 1.8 p. 198 in \cite{Stenstrom}). Since $I\subseteq \ker f,$ $f$ factors to a map $R/I\rightarrow M_{\ef_2}.$ But this map has to be zero as $R/I$ is torsion and $M_{\ef_2}$ is torsion-free in $\tau_2.$  Thus, $f$ is zero and so $d_2 (q_{12}(q))=q_{12}(d_1(q))$ for every $q\in M_{\ef_1}.$
\end{proof}

This lemma give us crucial ingredients of the proof of the following proposition.

\begin{prop} Suppose that a Gabriel filter $\ef_1$ is contained in a Gabriel filter $\ef_2$, that $\ef_1$ is differential and that $M$ is a right $R$-module. If either
\begin{itemize}
\item[i)] $M_{\ef_1}$ is torsion-free with respect to the torsion theory corresponding to $\ef_2,$  or

\item[ii)] $\ef_2$ is differential,
\end{itemize}
then any derivation on $M$ extends both to $M_{\ef_1}$ and $M_{\ef_2}$ in such a way that the extensions on $M_{\ef_1}$ and $M_{\ef_2}$ agree.
\label{Extending_to_right_modules_of_quotients}
\end{prop}

\begin{proof} In both cases, $d$ extends from $M$ to a derivation $d_1$ on $M_{\ef_1}$ by the differentiability of $\ef_1$ (part 2 of Proposition \ref{Golan_Proposition}).

In case i), $M_{\ef_1}$ is $\tau_2$-torsion-free, and so the kernel of the map $q_2: M_{\ef_1}\rightarrow (M_{\ef_1})_{\ef_2}$ is zero. Thus, we obtain the extension $\overline{d_{12}}$ of $d_1$ on $(M_{\ef_1})_{\ef_2}$ such that $q_2d_1=\overline{d_{12}}q_2.$ Since the map $i: (M_{\ef_1})_{\ef_2}\rightarrow M_{\ef_2}$ with $L_s\mapsto q_{12}(s)$ for  $s\in M_{\ef_1}$ is an isomorphism (see the diagram preceding Definition \ref{def_of_agreement}), the map $d_2:= i \overline{d_{12}}i^{-1}$ defines a derivation on $M_{\ef_2}$ such that
\[d_2q_{12}(s)=d_2i(L_s)=i\overline{d_{12}}(L_s)= i\overline{d_{12}}q_2(s)= i q_2 d_1(s)= i(L_{d_1(s)})=q_{12}d_1(s).\]
Thus the extensions agree by part i) of Lemma \ref{lemma_for_extending_right}.

In case ii), $d$ extends from $M$ to $M_{\ef_2}$ by differentiability of $\ef_2.$ Then the extensions on $M_{\ef_1}$ and $M_{\ef_2}$ agree by part ii) of Lemma \ref{lemma_for_extending_right}.
\end{proof}

\begin{lemma}
The relation of agreement of extensions is transitive. Namely, if $\ef_1\subseteq\ef_2\subseteq\ef_3$ are three Gabriel filters and $M$ a right $R$-module such that the extension on $M_{\ef_1}$ and $M_{\ef_2}$ agree and the extensions on
$M_{\ef_2}$ and $M_{\ef_3}$ agree, then the extensions on $M_{\ef_1}$ and $M_{\ef_3}$ agree also.
\label{transitivity}
\end{lemma}
\begin{proof}
Keeping the notation used in Lemma \ref{lemma_for_extending_right}, we have that $q_{13}=q_{23}q_{12},$ $d_2q_{12}=q_{12} d_1$ and $d_3q_{23}=q_{23} d_2.$ Then $d_3q_{13}=q_{13} d_1$ since
\[d_3q_{13}=d_3 q_{23}q_{12}= q_{23}d_2 q_{12}= q_{23}q_{12} d_1 =q_{13} d_1.\]
\end{proof}

The following corollary of Proposition \ref{Extending_to_right_modules_of_quotients} answers the questions Q1 and Q2 from the beginning of the section.
\begin{corollary}
Suppose that a Gabriel filter $\ef_1$ is contained in a Gabriel filter $\ef_2$. Let $Q_1$ and $Q_2$ denote the respective right rings of quotients. If either
\begin{itemize}
\item[i)] $\ef_1$ is differential and $Q_1$ is $\tau_2$-torsion-free,

\item[ii)] $\ef_1$ and $\ef_2$ are differential, or

\item[iii)] $\ef_1$ and $\ef_2$ are faithful,
\end{itemize}
then any derivation on $R$ extends both to $Q_1$ and $Q_2$ and the extensions agree.
\label{Extending_to_right_rings_of_quotients}
\end{corollary}
\begin{proof}
In cases i) and ii), the claim follows directly by Proposition \ref{Extending_to_right_modules_of_quotients}.

In case iii), any derivation on $R$ extends to both $Q_1$ and $Q_2$ by part 1 of Proposition \ref{Golan_Proposition}. Then the claim follows from part ii) of Lemma \ref{lemma_for_extending_right}.
\end{proof}

This corollary guarantees the agreement of extensions on some frequently considered right rings of quotients. One of them is the maximal right ring of quotients $\Qmax(R).$ This is the ring of quotients with respect to the Lambek torsion theory -- the torsion theory cogenerated by the injective envelope $E(R)$ (see sections 13B and 13C in \cite{Lam} and Example 1, page 200, \cite{Stenstrom}). The Lambek torsion theory is the largest faithful hereditary torsion theory. The Gabriel filter of this torsion theory is the set of all dense right ideals (see definition 8.2. in \cite{Lam} and Proposition VI 5.5, p. 147 in \cite{Stenstrom}). By Proposition 9 in \cite{Lia_Diff}, the Lambek torsion theory is differential.

If $R$ is a right Ore ring with the set of regular elements $T$ (i.e.,
$rT \cap tR \neq 0,$ for every $t \in T$ and $r\in R$), we can
define a hereditary torsion theory by the condition that a right
$R$-module $M$ is a torsion module iff for every $m\in M$, there
is a nonzero $t\in T$ such that $mt =0.$ This torsion theory is
called the classical torsion theory for a right Ore ring. It
is hereditary and faithful. Its right ring of quotients is called the classical right ring of quotients and is denoted by $\Qcl(R).$

Another frequently considered ring of quotients is related to the concept of perfect Gabriel filters. If a Gabriel filter $\ef$ has the property \[M_{\ef}\cong M\otimes_R R_{\ef}\] for every right $R$-module $M,$ then $\ef$ is called perfect and the right ring of quotients $R_{\ef}$ is called the perfect right ring of quotients. These filters are interesting as all modules of quotients are determined solely by the right ring of quotients. The classical torsion theory has this property. In fact, the perfect torsion theories have developed as generalization of the classical torsion theory. Perfect right rings of quotients and perfect filters have been studied and characterized in more details (see Theorem 2.1, page 227 in \cite{Stenstrom} and Proposition 3.4, page 231 in \cite{Stenstrom}). By Proposition 4 in \cite{Lia_Diff}, every perfect filter is differential.

Every ring has a maximal perfect right ring of quotients (see Theorem XI 4.1, p. 233, \cite{Stenstrom}). It is called the total right ring of quotients and is denoted by $\Qtot(R).$ $\Qtot(R)$ can be constructed as a directed union of all perfect right rings of quotients that are contained in $\Qmax(R)$ (see Theorem XI 4.1, p. 233, \cite{Stenstrom}). 
Note that $\Qtot(R)$ is the ring of quotients with respect to the torsion theory determined by $M\in \te$ iff $M\otimes_R \Qtot(R)=0$ (see Theorem 2.1, p. 227 and Proposition 3.4, p. 231 in \cite{Stenstrom}).

The class of nonsingular modules over a ring $R$ is closed
under submodules, extensions, products and injective envelopes so it is a torsion-free class of a hereditary torsion theory. This torsion theory is called the Goldie torsion theory. It
is larger than any hereditary faithful torsion theory (see Example
3, p. 26 in \cite{Bland_book}). So, the Lambek torsion theory is
smaller than the Goldie's. If $R$ is right nonsingular, the Lambek and Goldie torsion theories
coincide (see  p. 26 in \cite{Bland_book} or p. 149 in \cite{Stenstrom}). The Goldie torsion theory is differential (Proposition 14 in \cite{Lia_Diff}). Let us denote the ring of quotients with respect to the Goldie torsion theory by $\Qgold(R).$ It is isomorphic to the injective envelope of the torsion-free quotient of $R$  (Propositions IX 1.7, 2.5, 2.7, and 2.11 and Lemma IX 2.10 in \cite{Stenstrom}).

\begin{corollary} Let $\delta$ be a derivation on $R,$ $M$ a right $R$-module and $d$ a $\delta$-derivation on $M.$
\begin{enumerate}
\item The extension of $\delta$ on any right ring of quotients with respect to a hereditary and faithful torsion theory agrees with the extension of $\delta$ on $\Qmax(R).$ In particular, the extensions on $\Qtot(R)$ and $\Qmax(R)$ agree.

\item The extension of $d$ on any module of quotients of $M$ with respect to a differential, hereditary and faithful torsion theory agrees with the extension of $d$ on the module of quotients with respect to the Lambek torsion theory. In particular, the extensions on module of quotients of $M$ with respect to Lambek and maximal perfect torsion theory agree.

\item If $R$ is right Ore and $\Qcl(R)$ is its classical right ring of quotients, the extension on $\Qcl(R)$ agrees with the extensions on $\Qmax(R)$ and $\Qtot(R).$

\item The extensions of $\delta$ on $\Qmax(R)$ and $\Qgold(R)$ agree. Moreover, the extension of $\delta$ on any right ring of quotients with respect to a hereditary and faithful torsion theory agrees with the extension of $\delta$ on $\Qgold(R).$
\end{enumerate}
\label{extending_on_Qtot_and_Qmax}
\end{corollary}
\begin{proof}
(1) follows directly from Corollary \ref{Extending_to_right_rings_of_quotients} as the condition iii) of Corollary \ref{Extending_to_right_rings_of_quotients} is satisfied. Note that for $\Qtot(R)$ and $\Qmax(R)$ conditions i) and ii) are satisfied also (by Propositions 4 and 9 from \cite{Lia_Diff}).

(2) follows directly from part ii) of Lemma \ref{lemma_for_extending_right} using that the Lambek torsion theory is differential (Proposition 9 of \cite{Lia_Diff}). The second sentence of (2) holds since a perfect torsion theory is differential (Proposition 4 of \cite{Lia_Diff}).

(3) follows from (1) and Corollary \ref{Extending_to_right_rings_of_quotients} because the classical torsion theory is hereditary and faithful.

The first sentence of (4) holds by part ii) of Corollary \ref{Extending_to_right_rings_of_quotients}. Note that Lambek and Goldie torsion theories are differential by Propositions 9 and 14 of \cite{Lia_Diff}. The second sentence of (4) holds by (1) and transitivity of extensions (Lemma \ref{transitivity}).
\end{proof}

\section{Differentiability of symmetric torsion theories}
\label{section_differentiability_of_symmetric}

In this section, we shall work with derivations on bimodules. We shall prove the symmetric version of Golan's and Bland's results on differential one-sided torsion theories and the symmetric versions of results from \cite{Lia_Diff}.

If $R$ and $S$ are two rings, $\ef_l$ a Gabriel filter of left $R$-ideals and $\ef_r$ a Gabriel filters of right $S$-ideals, define $_{\ef_l}\ef_{\ef_r}$ as the set of right ideals of $S\otimes_{\Zset}
R^{op}$ containing an ideal of the form $J\otimes R^{op}+S\otimes I$ where $I\in\ef_l$ and $J\in \ef_r.$ This defines a Gabriel filter (\cite{Ortega_thesis}, page 100). We shorten the notation $_{\ef_l}\ef_{\ef_r}$ to
$_{l}\ef_{r}$ when there is no confusion about the Gabriel filters used. We call $_{l}\ef_{r}$ the {\em symmetric filter induced by $\ef_l$ and $\ef_r$}. If $\tau_l$ and $\tau_r$ are the torsion theories corresponding to filters $\ef_l$ and $\ef_r$ respectively, we call the torsion theory $_l\tau_r$ corresponding to $_{l}\ef_{r}$ {\em the torsion theory induced by $\tau_l$ and $\tau_r.$}

If $M$ is an $R$-$S$-bimodule, $\te_l(M),$ $\te_r(M)$ and $_l\te_r(M)$ torsion submodules of $M$ for $\tau_l,$ $\tau_r$ and  $_l\tau_r$ respectively, then \[_l\te_r(M)=\te_l(M)\cap\te_r(M).\] For details see \cite{Golan}, I, ch. 2, Proposition 2.5. Thus the torsion theory on $S\otimes_{\Zset}
R^{op}$ corresponding to filter $_l\ef_r$ of right $S\otimes_{\Zset}
R^{op}$-ideals is exactly the torsion theory of $R$-$S$-bimodules with the torsion class $\te_l\cap\te_r.$ Thus the following lemma holds.

\begin{lemma} Using the notation as above, the following conditions are equivalent for every $x\in M.$
\begin{enumerate}
\item $x\in\, _l\te_r(M)$

\item $\ann_l(x)\in \ef_l$ and $\ann_r(x)\in \ef_r.$

\item $\ann_r(x)\otimes_{\Zset} R^{op}+R\otimes_{\Zset}\ann_l(x)\in\,_l\ef_r.$

\item $_l\ann_r(x):=\{t\in R\otimes R^{op}| xt=0\}\in\,_l\ef_r.$
\end{enumerate}
\label{torsion_vs_filter}
\end{lemma}
\begin{proof}
By definition of induced filters and torsion theory and since $_l\te_r(M)=\te_l(M)\cap\te_r(M)$ we have that $(1)\Rightarrow (2)\Rightarrow (3)\Rightarrow (4).$ (4) states that the right annihilator of $x$ in $S\otimes_{\Zset}R^{op}$ is in $_l\ef_r$ and this implies that $x$ is in torsion submodule with respect to the torsion theory $_l\tau_r.$ So we have (1) in this case.
\end{proof}

In \cite{Ortega_paper} and \cite{Ortega_thesis}, Ortega defines the {\em symmetric module of quotients} $_{\ef_l}M_{\ef_r}$ of $M$ with respect to $_l\ef_r$ to be
\[_{\ef_l}M_{\ef_r}=\dirlim_{K\in _l\ef_r}\;  \homo(K, \frac{M}{_l\te_r(M)})\]
where the homomorphisms in the formula are $S\otimes R^{op}$ homomorphisms (equivalently $R$-$S$-bimodule homomorphisms).
We shorten the notation $_{\ef_l}M_{\ef_r}$ to $_lM_r.$ Just as in the right-sided case, there is a left exact functor $q_M$ mapping $M$ to the symmetric module of quotients $_lM_r$ such that $\ker q_M$ is the torsion module $_l\te_r(M)$ (see Lemma 3.1 in \cite{Lia_Sym}).

If $R=S,$ the module $_lR_r$ has a ring structure (see \cite{Schelter} or Lemma 1.5 in \cite{Ortega_paper}). The ring $_lR_r$ is called {\em the symmetric ring of quotients} with respect to the torsion theory $_l\tau_r.$

\begin{example}
\begin{enumerate}
\item Let $\de_r$ and $\de_l$ denote the Gabriel filters of all dense right and left $R$-ideals respectively.  Let $_l\de_r$ denote the symmetric filter induced by $\de_r$ and $\de_l$. The corresponding symmetric ring of quotients is called the maximal symmetric ring of quotients. We denote it by $\Qsimmax(R).$ Utumi first studied this ring in \cite{Utumi}. Lanning (\cite{Lanning}), Schelter (\cite{Schelter}) and Ortega (\cite{Ortega_paper_Qsimmax}, \cite{Ortega_thesis}) studied $\Qsimmax(R)$ using torsion theories. $\Qsimmax(R)$ can also be described as follows:
\[\begin{array}{rcl}\Qsimmax(R) & = & \{q\in \Qmax(R)| Iq\subseteq R\mbox{ for some }I\in\de_l\}\\
& = & \{q\in \Qlmax(R)| qJ\subseteq R\mbox{ for some }J\in\de_r\}.\end{array}\]
For proof see Remark 4.33 in \cite{Ortega_thesis}.

As $R$ is torsion-free for the torsion theories determined by $\de_r$ and $\de_l$, $R$ is torsion-free in the torsion theory induced by $\de_r$ and $\de_l$ as well. Thus, $R$ embeds in $\Qsimmax(R).$ $_{\ef_l}R_{\ef_r}$ also embeds in $\Qsimmax(R)$ for every $_l\tau_r$ induced by hereditary and faithful left and right torsion theories $\tau_l$ and $\tau_r.$

\item If $R$ is a right and left Ore ring with the set of regular elements $T$, the
classical left and the classical right rings of quotients coincide (see page 303 in \cite{Lam}) and we denote this ring of quotients by $\Qlrcl(R).$ The torsion theory for $R$-bimodules defined by the condition that an $R$-bimodule $M$ is a torsion module iff for every $m\in M$, there are nonzero $t,s\in T$ such that $sm=mt=0$ is induced by the left and right classical torsion theories. This torsion theory is called the classical torsion theory for an Ore ring. $\Qlrcl(R)$ is the symmetric ring of quotients with respect to the classical torsion theory (see Corollary 5.2 in \cite{Lia_Sym}).

\item The symmetric Goldie torsion theory is the theory induced by the torsion theories of left and right $R$-modules whose torsion-free classes consist of left and right nonsingular modules respectively. We denote the symmetric ring of quotients with respect to this torsion theory by $\Qsimgold(R).$

\item In \cite{Lia_Sym}, the symmetric version of a right perfect ring of quotients and the symmetric version of the total right ring of quotients are defined and studied. A ring homomorphism $f: R\rightarrow S$ makes $S$ into a perfect symmetric ring of quotients if the family of left ideals $\ef_l=\{I | Sf(I)=S\}$ is a
left Gabriel filter, the family of right ideals $\ef_r=\{J | f(J)S=S\}$ is a right
Gabriel filter and there is a ring isomorphism $g: S\cong\; _{\ef_l}R_{\ef_r}$ such that
$g\circ f$ is the canonical map $q_R: R\rightarrow\; _{\ef_l}R_{\ef_r}.$ For the equivalent versions of this definition see Theorem 4.1 of \cite{Lia_Sym}. In this case, the filter $_l\ef_r$ induced by $\ef_l$ and $\ef_r$ is called perfect and the torsion theory it determines is called a perfect symmetric torsion theory.

Symmetrizing the condition for the perfect right torsion theories $M_{\ef}\cong M\otimes_R R_{\ef}$, we arrive to the condition describing the perfect symmetric torsion theories and symmetric filters:
\[_lR_r \otimes_R M\otimes_R\, _lR_r\cong \;_lM_r\]
for all $R$-bimodules $M.$ For the proof of this fact together with equivalent conditions for describing perfect symmetric torsion theories and filters see Theorem 4.2 in \cite{Lia_Sym}.

By Theorem 5.1 in \cite{Lia_Sym}, every ring has the largest perfect symmetric ring of quotients -- the total symmetric ring of quotients. We denote it by $\Qsimtot(R).$ Analogously to the right-sided case, $R\subseteq \Qsimtot(R)\subseteq \Qsimmax(R).$ If $R$ is Ore, $R\subseteq\Qlrcl(R)\subseteq \Qsimtot(R)\subseteq \Qsimmax(R)$ (see Corollary 5.2 in \cite{Lia_Sym}).
\end{enumerate}
\label{example_of_symmetric}
\end{example}

Let us turn to ring derivations now. Note that every derivation $\delta$ on $R$ determines a derivation on $R\otimes_{\Zset}R^{op}$ given by \[\overline{\delta}(r\otimes s)=\delta(r)\otimes s+r\otimes \delta(s).\]

If $M$ is an $R$-bimodule, and $\delta$ a derivation on $R$, we shall say that an additive map $d: M\rightarrow M$ is a {\em $\delta$-derivation} if
\[d(xr)=d(x)r+x\delta(r)\mbox{ and }d(rx)=\delta(r)x+rd(x)\]
for all $x\in M$ and $r\in R.$ Note that $d$ is a $\overline{\delta}$-derivation on $M$ considered as a right $R\otimes_{\Zset}R^{op}$-module since

\[\begin{array}{rcl}
d(x(r\otimes s)) & = & d(sxr)=\delta(s)xr+sd(x)r+sx\delta(r)\\
& = & x(r\otimes\delta(s))+d(x)(r\otimes s)+ x(\delta(r)\otimes s)\\
& = & d(x)(r\otimes s) + x\overline{\delta}(r\otimes s).
\end{array}\]

Conversely, every $\overline{\delta}$-derivation of a right $R\otimes_{\Zset}R^{op}$-module determines a $\delta$-de\-ri\-va\-tion of a bimodule:

\[\begin{array}{rcl}
d(xr) & = & d(x(r\otimes 1)=d(x)(r\otimes 1)+x\overline{\delta}(r\otimes1)\\
& = & 1d(x)r +x(\delta(r)\otimes 1)+ x(r\otimes\delta(1))\\
& = & d(x)r + 1 x\delta(r)+0 = d(x)r + x\delta(r). \end{array}\]
$d(rx)=\delta(r)x+rd(x)$ follows similarly.

Thus, every derivation $\delta$ on $R$ is a $\overline{\delta}$-derivation on $R$ considered as a right $R\otimes_{\Zset}R^{op}$-module. Conversely, every derivation $\overline{\delta}$ on $R\otimes_{\Zset}R^{op}$ is a $\delta$-derivation of $R\otimes_{\Zset}R^{op}$ considered as an $R$-bimodule.

Consider a symmetric filter $_l\ef_r$ induced by a left Gabriel filter $\ef_l$ and a right Gabriel filter $\ef_r.$ We shall say that $_l\ef_r$ is a {\em differential filter}
if for every $I\in\, _l\ef_r$ there is $J\in\, _l\ef_r$ such that $\overline{\delta}(J)\subseteq I$ for all $R\otimes_{\Zset}R^{op}$ derivations $\overline{\delta}.$ If we consider the right $R\otimes_{\Zset}R^{op}$-ideals $I$ and $J$ as $R$-bimodules, the condition $\overline{\delta}(J)\subseteq I$ is equivalent with $\delta(J)\subseteq I$ by observations above. The hereditary torsion theory determined by $_l\ef_r$ is said to be {\em differential} in this case.

The equivalence of the three conditions in the proposition below produces the symmetric version of Lemma 1.5 from \cite{Bland_paper}. In what follows, we also use the notation $_l\ann_r(x)$ for $\{t\in R\otimes R^{op}| xt=0\}.$
\begin{prop}
If $_l\ef_r$ is a Gabriel filter induced by $\ef_l$ and $\ef_r$ corresponding to the torsion theory $_l\tau_r,$ the following conditions are equivalent.
\begin{enumerate}
\item $_l\ef_r$ is a differential filter.

\item For every $R$-bimodule $M,$ every ring derivation $\delta$ and $\delta$-derivation $d$ on $M,$ $d(\,_l\te_r(M))\subseteq\, _l\te_r(M).$

\item For every $R$-bimodule $M,$ every ring derivation $\delta$ and $x\in\,_l\te_r(M),$ there is $K\in\,_l\ef_r$ such that $\overline{\delta}(K)\subseteq\, _l\ann_r(x).$

Moreover, if  $\ef_l$ and $\ef_r$ are differential, then $_l\ef_r$ is also differential.
\end{enumerate}
\label{symmetric_differential_TT}
\end{prop}
\begin{proof}
$(1)\Rightarrow (2)$ Let $x\in\, _l\te_r(M).$ By Lemma \ref{torsion_vs_filter} $_l\ann_r(x)$ is in $_l\ef_r.$ By assumption, there is $K\in\,_l\ef_r$ such that $\overline{\delta}(K)\subseteq\,_l\ann_r(x).$ Let us consider the ideal $I=K\cap\,_l\ann_r(x).$ By construction, every $t\in I$ is such that $xt=0$ and $x\overline{\delta}(t)=0.$ Thus $d(xt)=0$ and so $d(x)t=d(xt)-x\overline{\delta}(t)=0.$ So, $I\subseteq\,_l\ann_r(d(x)).$ As $I\in\,_l\ef_r,$ $_l\ann_r(d(x))\in\,_l\ef_r.$ But then $d(x)$ is in  $_l\te_r(M)$ by Lemma \ref{torsion_vs_filter}.

$(2)\Rightarrow (3)$ Let $x\in\, _l\te_r(M).$ By assumption $d(x)$ is in $_l\te_r(M)$ then too. By
Lemma \ref{torsion_vs_filter} $\ann_l(x), \ann_l(d(x))\in \ef_l$ and $\ann_r(x), \ann_r(d(x))\in \ef_r$ so
$I=\ann_l(x)\cap \ann_l(d(x))\in \ef_l$ and $J=\ann_r(x)\cap\ann_r(d(x))\in \ef_r.$ Then $K=J\otimes_{\Zset} R^{op}+R\otimes_{\Zset}I$ is in $_l\ef_r.$

Note that a generator $r\otimes s$ of
$J\otimes_{\Zset} R^{op}$ is such that $x(r\otimes s)=sxr=s0=0$ and $d(x)(r\otimes s)=sd(x)r=s0=0.$ Similarly, a generator $r\otimes s$ of $R\otimes_{\Zset}I$ is in right annihilators of $x$ and $d(x)$ in $R\otimes_{\Zset}R^{op}.$
Thus, an element $t$ from $K$ is such that $xt=0$ and $d(x)t=0.$ Then $x\overline{\delta}(t)=d(xt)-d(x)t=0.$ So, $\overline{\delta}(K)\subseteq\, _l\ann_r(x).$

$(3)\Rightarrow (1)$ If $I\in\,_l\ef_r,$ then $(1\otimes 1)+I$ is in $_l\te_r((R\otimes_{\Zset}R^{op})/I).$ By assumption, there is $J\in\,_l\ef_r$ such that $\overline{\delta}(J)\subseteq\, _l\ann_r((1\otimes 1)+I).$ As $_l\ann_r((1\otimes 1)+I)=I,$ $\overline{\delta}(J)\subseteq I.$

Let us now assume that $\ef_l$ and $\ef_r$ are differential. Let $M$ be an $R$-bimodule, $\delta$ a ring derivation and $d$ a $\delta$-derivation on $M.$ As $\ef_l$ and $\ef_r$ are differential,
$d(\te_l(M))\subseteq \te_l(M)$ and $d(\te_r(M))\subseteq \te_r(M)$ by Proposition \ref{Golan_Proposition}. Thus
$d(_l\te_r(M))=d(\te_l(M)\cap \te_r(M)))\subseteq d(\te_l(M))\cap d(\te_r(M))$ is contained in $\te_l(M)\cap\te_r(M)=_l\te_r(M)$ and so $_l\ef_r$ is differential.
\end{proof}

The following theorem is the symmetric version of Bland's and Golan's results (Proposition \ref{Golan_Proposition} and Theorem \ref{Bland_Theorem}).

\begin{theorem} Let $\delta$ be a derivation on $R$, $M$ an $R$-bimodule, $d$ a $\delta$-derivation on $M,$ and $_l\ef_r$ a Gabriel filter induced by $\ef_l$ and $\ef_r$ corresponding to the torsion theory $_l\tau_r.$
\begin{enumerate}
\item If $M$ is torsion-free, then $d$ extends to a derivation on the module of quotients $_lM_r$ such that $d q_M=q_M d.$

\item If $d(\, _l\te_r(M))\subseteq\, _l\te_r(M),$ then $d$ extends to a derivation on the module of quotients $_lM_r$ such that $d q_M=q_M d.$

\item If a derivation on a module $M$ extends to a derivation on the module of quotients $_lM_r,$ then such extension is unique.

\item The filter $_l\ef_r$ is differential if and only if
every derivation on any module $M$ extends uniquely to a derivation on the module of quotients $_lM_r.$
\end{enumerate}
\label{extending_symmetric_diff}
\end{theorem}
\begin{proof}
(1) Consider $M$ as a right $R\otimes_{\Zset}R^{op}$-module. Note that $_lM_r$ is a right $R\otimes_{\Zset}R^{op}$-module of quotients with respect to $_l\tau_r$ considered as the torsion theory of right $R\otimes_{\Zset}R^{op}$-modules. Thus,
$d$ extends to a $\overline{\delta}$-derivation on $_lM_r$ with $d q_M=q_M d$ by condition (1) of Proposition \ref{Golan_Proposition}. But then such extension is a $\delta$-derivation on $_lM_r.$

(2) As $d(\, _l\te_r(M))\subseteq\, _l\te_r(M),$ $d$ defines a derivation on $M/\,_l\te_r(M).$ As $_l\f_r(M)=M/\,_l\te_r(M)$ is torsion-free, this derivation extends to the quotient $_l(\,_l\f_r(M))_r$ by $d(m+\,_l\te_r(M))=d(m)+\,_l\te_r(M).$ But $_l(\,_l\f_r(M))_r=\,_lM_r$ (see part 1 of Lemma 3.1 in \cite{Lia_Sym}), so we obtain the extension of $d$ on $_lM_r$ with
$d q_{\f(M)}=q_{\f(M)} d.$ As $q_M(m)=q_{\f(M)}(m+\,_l\te_r(M)),$
\[\begin{array}{rcl}
d (q_M(m)) & = & d(q_{\f(M)}(m+\,_l\te_r(M)))\\
& = & q_{\f(M)}(d((m+\,_l\te_r(M))))\\
& = & q_{\f(M)}(d(m)+\,_l\te_r(M))\\
& = & q_M(d(m)).
\end{array}\]
Thus $d q_M=q_M d.$

(3) Let us assume that we have two extensions $d_1$ and $d_2$ on $_lM_r$ extending $d$ on $M.$ As $d_1-d_2$ is an $R$-bimodule homomorphism with $M\subseteq\ker(d_1-d_2),$ $d_1-d_2$ factors to a homomorphism of $_lM_r/M\rightarrow\,_lM_r.$ However, $_lM_r/M$ is torsion and $_lM_r$ torsion-free, so $d_1-d_2$ is zero on $_lM_r/M.$ Hence, $d_1=d_2$ on $_lM_r.$

(4) If $_l\ef_r$ is differential, it follows from (2), (3) and Proposition \ref{symmetric_differential_TT} that a derivation on a module extends uniquely to a derivation on the module of quotients. Conversely, if a derivation on a bimodule $M$ extends to $_lM_r,$ then every $m\in\,_l\te_r(M)$ maps to zero by $d q_M=q_M d$ as $q_M(m)=0$ Thus, $q_M(d(x))=0$ and so $d(x)$ is in $\ker q_M=\,_l\te_r(M).$ Hence $d(\,_l\te_r(M))\subseteq\,_l\te_r(M)$ and $_l\ef_r$ is differential by Proposition \ref{symmetric_differential_TT}.
\end{proof}

In \cite{Lia_Diff}, it is shown that one-sided Lambek and Goldie torsion theories are differential. This fact has the following corollary.

\begin{corollary}
The symmetric Lambek and Goldie torsion theories are differential.
\label{symmetric_Lambek_Goldie}
\end{corollary}
\begin{proof}
As Lambek and Goldie torsion theories of right (and left) modules are differential by Propositions 9 and 14 in \cite{Lia_Diff} and symmetric Lambek and Goldie torsion theories are induced by one-sided Lambek and Goldie theories respectively (see Example \ref{example_of_symmetric}), the result directly follows by the last part of Proposition \ref{symmetric_differential_TT}.
\end{proof}

Another result from \cite{Lia_Diff} is that every one-sided perfect torsion theory is differential. We prove the symmetric version of this result.
\begin{corollary}
If a symmetric filter is perfect, it is differential.
\label{symmetric_perfect}
\end{corollary}
\begin{proof}
Recall that a ring $S$ with a ring homomorphism $f: R\rightarrow S$ is a perfect symmetric ring of quotients if the family of left ideals $\ef_l=\{I | Sf(I)=S\}$ is a
left Gabriel filter, the family of right ideals $\ef_r=\{J | f(J)S=S\}$ is a right
Gabriel filter and there is a ring isomorphism $g: S\cong\; _{\ef_l}R_{\ef_r}$ such that
$g\circ f$ is the canonical map $q_R: R\rightarrow\; _{\ef_l}R_{\ef_r}$ (see Theorem 4.1 of \cite{Lia_Sym}). In this case, $\ef_r$ is a perfect right, $\ef_l$ a perfect left and $_l\ef_r$ a perfect symmetric filter. By Proposition 4 of \cite{Lia_Diff}, $\ef_l$ and $\ef_r$ are differential. Then $_l\ef_r$ is differential by the last part of Proposition \ref{symmetric_differential_TT}.
\end{proof}

As a consequence, every derivation on an $R$-bimodule $M$ lifts uniquely to a derivation on the module of quotients $_lM_r$ if $_l\ef_r$ is a Gabriel filter of the Lambek, Goldie or a perfect symmetric torsion theory. In particular, every derivation on $R$ uniquely lifts to a derivation of $\Qsimmax(R),$ $\Qsimtot(R)$ and $\Qsimgold(R).$

\section{Extending derivations to symmetric rings and modules of quotients}
\label{section_extending_on_symmetric}

In this section, we study the agreement of extensions of derivations on symmetric rings and modules of quotients. We prove the symmetric versions of results from section \ref{section_extending_on_right}.

Let $\ef^1_l$ and $\ef^2_l$ be left and $\ef^1_r$ and $\ef^2_r$ be right Gabriel filters inducing the symmetric filters $_l\ef_r^1$ and $_l\ef_r^2.$ Suppose that $_l\ef_r^1$ is contained in $_l\ef_r^2$ and that $M$ is an $R$-bimodule. Let $_l\tau_r^i$ denote the corresponding torsion theories and $q_i$ denote the left exact functors mapping $M$ to the rings of quotients $_lM_r^i$ with respect to $_l\ef_r^i$ for $i=1$ and $2.$ Note that we have the mapping $q_{12}:\,_lM_r^1\rightarrow\,_lM_r^2$ induced by the inclusion $_l\ef_r^1\subseteq\,_l\ef_r^2$ such that $q_{12}q_1=q_2$ just as in the right-sided case.

\begin{definition}
Suppose that there is a $\delta$-derivation $d$ defined on an $R$-bimodule $M.$ If $d$ extends to derivations $d_1$ on $_lM_r^1$ and $d_2$ on  $_lM_r^2$ such that the diagram
\begin{diagram}
  &             & _lM_r^1      &     & \rTo^{d_1} &             & _lM_r^1\\
  & \ruTo^{q_1} & \vLine       &     &            & \ruTo^{q_1} &    \\
M &             & \rTo^{d}     &     & M          &             & \dTo_{q_{12}}\\
  & \rdTo^{q_2} & \dTo_{q_{12}}&     &            & \rdTo^{q_2} &     \\
  &             & _lM_r^2      &     & \rTo^{d_2} &             & _lM_r^2\\
\end{diagram}
commutes, we say that the extensions of $d$ on $_lM_r^1$ and $_lM_r^2$ {\em agree.}
\end{definition}

The following are the symmetric versions of Lemma \ref{lemma_for_extending_right}, Proposition \ref{Extending_to_right_modules_of_quotients} and Corollary \ref{Extending_to_right_rings_of_quotients}.

\begin{lemma} Suppose that a symmetric filter $_l\ef_r^1$ is contained in a symmetric filter $_l\ef_r^2$ and that $M$ is an $R$-bimodule with a $\delta$-derivation $d$. If $d$ can be extended to $_lM_r^1$ and either
\begin{itemize}
\item[i)] $d$ can be extended from $_lM_r^1$ to $_lM_r^2,$  or

\item[ii)] $d$ can be extended from $M$ to $_lM_r^2,$
\end{itemize}
then the extensions of $d$ to $_lM_r^1$ and $_lM_r^2$ agree.
\label{lemma_for_extending_symmetric}
\end{lemma}
\begin{proof}
The proof follows the proof of Lemma \ref{lemma_for_extending_right} exactly using Theorem \ref{extending_symmetric_diff} instead of Theorem \ref{Bland_Theorem}.
\end{proof}

\begin{prop}
Suppose that a symmetric Gabriel filter $_l\ef_r^1$ is contained in a symmetric Gabriel filter $_l\ef_r^2$, that $_l\ef_r^1$ is differential and that $M$ is a right $R$-module. If either
\begin{itemize}
\item[i)] $_lM_r^1$ is torsion-free with respect to the torsion theory corresponding to $_l\ef_r^2,$  or

\item[ii)] $_l\ef_r^2$ is differential,
\end{itemize}
then any derivation on $M$ extends both to $_lM_r^1$ and $_lM_r^2$ in such a way that the extensions on $_lM_r^1$ and $_lM_r^2$ agree.
For $M=R,$ i), ii) or
\begin{itemize}
\item[iii)] $_l\ef_r^1$ and $_l\ef_r^2$ are faithful,
\end{itemize}
imply that any derivation on $R$ extends both to $_lR_r^1$ and $_lR_r^2$ in such a way that the extensions on $_lR_r^1$ and $_lR_r^2$ agree.
\label{extending_on_symmetric_modules}
\end{prop}
\begin{proof}
The proof follows exactly the proofs of Proposition \ref{Extending_to_right_modules_of_quotients} and Corollary \ref{Extending_to_right_rings_of_quotients} using Proposition \ref{symmetric_differential_TT}, Theorem \ref{extending_symmetric_diff} and Lemma \ref{lemma_for_extending_symmetric} in places of Lemma 1.5 from \cite{Bland_paper}, Proposition \ref{Golan_Proposition},  Theorem \ref{Bland_Theorem} and Lemma \ref{lemma_for_extending_right}.
\end{proof}

Note also that the relation of agreement is transitive (i.e. the symmetric version of Lemma \ref{transitivity} holds).
The proof in symmetric case is analogous to the proof of the right-sided case.

In Example \ref{example_of_symmetric}, we have recalled the definitions of the maximal symmetric ring of quotients $\Qsimmax(R)$, total symmetric ring of quotients $\Qsimtot(R),$ symmetric ring of quotients with respect to the symmetric Goldie torsion theory $\Qsimgold(R)$ and the concept of perfect symmetric ring of quotients and perfect symmetric filter. Now we prove the symmetric version of Corollary \ref{extending_on_Qtot_and_Qmax}.

\begin{corollary}
Let $\delta$ be a derivation on $R,$ $M$ an $R$-bimodule and $d$ a $\delta$-derivation on $M.$
\begin{enumerate}
\item The extension of $\delta$ on any right ring of quotients with respect to a hereditary and faithful symmetric torsion theory agrees with the extension of $\delta$ on $\Qsimmax(R).$ In particular, the extensions on $\Qsimtot(R)$ and $\Qsimmax(R)$ agree.

\item The extension of $d$ on any module of quotients of $M$ with respect to a differential, hereditary and faithful symmetric torsion theory agrees with the extension of $d$ on the module of quotients with respect to the symmetric Lambek torsion theory. In particular, the extensions on module of quotients of $M$ with respect to symmetric Lambek and maximal perfect symmetric torsion theory agree.

\item If $R$ is Ore and $\Qlrcl(R)$ is its classical ring of quotients, the extension on $\Qlrcl(R)$ agrees with the extensions on $\Qsimmax(R)$ and $\Qsimtot(R).$

\item The extensions of $\delta$ on $\Qsimmax(R)$ and $\Qsimgold(R)$ agree. Moreover, the extension of $\delta$ on any right ring of quotients with respect to a hereditary and faithful torsion theory agrees with the extension of $\delta$ on $\Qsimgold(R).$
\end{enumerate}
\label{extending_on_symmetric_Lambek_Goldie_perfect}
\end{corollary}
\begin{proof}
(1) follows directly from Proposition \ref{extending_on_symmetric_modules} as the condition iii) of Proposition \ref{extending_on_symmetric_modules} is satisfied. Note that for $\Qsimtot(R)$ and $\Qsimmax(R)$ conditions i) and ii) are satisfied also (by Corollaries \ref{symmetric_Lambek_Goldie} and \ref{symmetric_perfect}).

(2) follows directly from part ii) of Lemma \ref{lemma_for_extending_symmetric} using that the symmetric Lambek torsion theory is differential (Corollary \ref{symmetric_Lambek_Goldie}). The second sentence holds by ii) of Proposition \ref{extending_on_symmetric_modules} since a perfect symmetric torsion theory is differential (Corollary \ref{symmetric_perfect}).

(3) follows from (1) and Proposition \ref{extending_on_symmetric_modules} because the classical symmetric torsion theory is hereditary and faithful (and also perfect).

The first sentence of (4) holds by Proposition \ref{extending_on_symmetric_modules} since Lambek and Goldie torsion theories are differential by Corollary \ref{symmetric_Lambek_Goldie}. The second sentence of (4) holds by (1) and transitivity of extensions.
\end{proof}

\end{document}